\documentclass[11pt]{article}

\usepackage[french]{babel}
\usepackage[T1]{fontenc}
\usepackage{color}

\newtheorem{defi}{D\'efinition}
\newtheorem{prop}[defi]{Proposition}
\newtheorem{theo}[defi]{Th\'eor\`eme}
\newtheorem{conj}[defi]{Conjecture}
\newtheorem{lemm}[defi]{Lemme}
\newtheorem{coro}[defi]{Corollaire}
\newtheorem{rema}[defi]{Remarque}
\newtheorem{exem}[defi]{Exemple}
\newtheorem{exer}{Exercice}%[section]

\newcommand{\bdefi}{\begin{defi}}
\newcommand{\edefi}{\end{defi}}
\newcommand{\bprop}{\begin{prop}}
\newcommand{\eprop}{\end{prop}}
\newcommand{\btheo}{\begin{theo}}
\newcommand{\etheo}{\end{theo}}
\newcommand{\blemm}{\begin{lemm}}
\newcommand{\brema}{\begin{rema}}
\newcommand{\erema}{\end{rema}}
\newcommand{\bexer}{\begin{exer}}
\newcommand{\eexer}{\end{exer}}
\newcommand{\bconj}{\begin{conj}}
\newcommand{\econj}{\end{conj}}
\newcommand{\elemm}{\end{lemm}}
\newcommand{\bcoro}{\begin{coro}}
\newcommand{\ecoro}{\end{coro}}
\newcommand{\bexem}{\begin{exem}}
\newcommand{\eexem}{\end{exem}}
\newcommand{\dem}{\noindent{\bf Preuve. }}

\newcommand{\R}{{\cal R}}

\newcommand{\M}{{\cal M}}
\renewcommand{\P}{{\cal P}}

\newcommand{\E}{{\cal E}}

\newcommand{\C}{{\cal C}}

\newcommand{\Q}{{\cal Q}}
\newcommand{\F}{{\cal F}}
\newcommand{\V}{{\cal V}}

\usepackage{amssymb}

\newcommand{\maths}[1]{{\mathbb #1}}  

\newcommand{\RR}{\maths{R}}
\newcommand{\NN}{\maths{N}}
\newcommand{\CC}{\maths{C}}
\newcommand{\QQ}{\maths{Q}}

\renewcommand{\SS}{\maths{S}}

\newcommand{\HH}{\maths{H}}

\newcommand{\ZZ}{\maths{Z}}

\newcommand{\PP}{\maths{P}}

\newcommand{\ra}{\rightarrow}

\newcommand{\ov}[1]{{\overline{#1}}} 
\newcommand{\wt}[1]{{\widetilde{#1}}}

\newcommand{\ga}{\gamma}
\newcommand{\Ga}{\Gamma}

\newcommand{\rfd}{\mbox{$\R_{fdp}(\pi_1S,{\rm PSL}_2(\RR))/%
{\rm PSL}_2(\RR)$}}
\newcommand{\wts}{\mbox{$\wt{\rm Teich}(S)$}}
\newcommand{\ts}{\mbox{${\rm Teich}(S)$}}
\newcommand{\tei}{Teich\-mül\-ler}
\newcommand{\cinf}{{{\rm C}^\infty}}
\newcommand{\cqfd}{\hfill$\Box$}

\title{Sur la compactification de Thurston \\
de l'espace de Teichmüller}

\author{Frédéric Paulin}

\date{}

\begin{document}
\maketitle
\begin{abstract}
The aim of these lecture notes is, after having quickly described
various compactifications of the Teichm\"uller space of a
compact connected oriented surface minus finitely many points, to give
a construction, by the equivariant Gromov topology introduced by the
author, of Thurston's compactification of this Teichm\"uller space.  
\footnote{ 
    Riemann surface; Thurston compactification; Teichmüller space;
    hyperbolic surface; equivariant Gromov topology.  }
\end{abstract} 

\bigskip
Le but de ces notes d'exposés est de donner une construction, par la
topologie de Gromov équivariante introduite dans \cite{Pau2}, de la
compactification de Thurston de l'espace de \tei\ $\ts$ d'une surface
$S$ orientée connexe compacte privée d'un nombre fini de points.

\bigskip
Il existe de nombreuses compactifications (géométriquement
intéressantes) de l'espace de \tei\ de $S$. Avant de donner en détail
notre approche de la compactification de Thurston, nous décrivons,
pour la culture et de manière non exhaustive, les principales
compactifications de $\ts$, de manière briève (et donc sans pouvoir
donner toutes les explications nécessaires, pour lesquelles nous
renvoyons à la bibliographie), ainsi que leurs relations, en supposant
pour simplifier que $S$ est compacte.

\medskip
\begin{itemize}
\item[$\bullet$] {\bf Les compactifications de \tei} (voir
  \cite{Abi,Gar,Nag}). Il y en a une pour toute structure de surface
  de Riemann $X$ fixée sur $S$, par l'espace des points à l'infini des
  rayons de \tei\ issus de (l'image encore notée $X$ dans $\ts$ de)
  $X$. La distance de \tei\ est finslérienne complète sur la variété
  différentielle $\ts$, et ces rayons de \tei\ sont les rayons
  géodésiques issus de $X$ pour cette métrique, ce qui fournit la
  compactification de \tei\ de $\ts$ associée à $X$. Le bord
  s'identifie avec l'espace $\Q_1(X)$ des formes différentielles
  quadratiques holomorphes de norme $1$ sur $X$, qui s'identifie à la
  sphère unité de l'espace cotangent de $\ts$ en $X$.

\item[$\bullet$] {\bf Les compactifications de Bers} (voir par exemple
  \cite{Ber,Nag,Bro,McM}). Il y en a une pour toute structure de
  surface de Riemann $X$ fixée sur $S$, construite comme suit. La
  surface de Riemann $X$ est isomorphe à la surface de Riemann
  quotient $\Ga_X\backslash\HH$ où $\HH\subset\CC$ est le demi-plan
  supérieur, et $\Ga_X$ un sous-groupe discret sans torsion du groupe
  ${\rm Aut}(\HH)\simeq {\rm PSL}_2(\RR)\subset{\rm PSL}_2(\CC) $ des
  automorphismes complexes de $\HH$.  Pour toute autre structure de
  surface de Riemann $Y$ sur $S$, soit $h:X\ra Y$ un homéomorphisme
  quasi-conforme homotope à l'identité envoyant la structure de $X$
  sur celle de $Y$.  Sa différentielle de Beltrami
  $\overline{\partial} h/\partial h$, qui est presque partout définie,
  se relève en une différentielle de Beltrami $\Ga_X$-invariante sur
  $\HH$, qui étendue par $0$ sur la sphère de Riemann $\widehat{\CC}$,
  donne une différentielle de Beltrami $\mu$ sur $\widehat{\CC}$, qui
  est $\Ga_X$-invariante. Le théorème d'Ahlfors-Bers dit qu'il existe
  un homéomorphisme quasi-conforme $f_\mu:
  \widehat{\CC}\ra\widehat{\CC}$, de différentielle de Beltrami $\mu$,
  uniquement défini modulo conjugaison au but par un élément de ${\rm
    Aut}(\widehat{\CC})$, qui conjugue $\Ga_X$ à un sous-groupe
  (encore discret et sans torsion) de ${\rm Aut}(\widehat{\CC})\simeq
  PSL_2(\CC)$.  Ainsi l'application $$Y\mapsto \{\ga\mapsto f_\mu\circ
  \ga\circ f_\mu^{-1}\}$$ induit un homéomorphisme sur son image, de
  l'espace de \tei\ de $S$, vu comme l'espace quotient de l'espace des
  structures de surfaces de Riemann sur $S$ modulo l'action de ${\rm
    Diff}_0(S)$, à valeurs dans l'espace topologique quotient (non
  séparé) ${\rm Hom}(\Ga_X,{\rm PSL}_2(\CC))/{\rm PSL}(\CC)$, dont
  l'image est d'adhérence (séparée et) compacte, ce qui fournit une
  compactification de $\ts$.

\item[$\bullet$] {\bf La compactification naturelle de Thurston} (voir
  \cite{Thu0, FLP} pour la construction originale, ainsi que les
  rappels ci-dessous). Le bord est l'espace $\P\M\F(S)$ des classes
  d'équivalence de feuilletages transversalement mesurés à
  singularités de type selle à $k\geq 3$ branches, modulo isotopies,
  opérations de Whitehead et multiplication par une constante
  strictement positive de la mesure transverse.  Une propriété
  cruciale de la compactification de Thurston est qu'elle est
  naturelle (au sens rappelé ci-dessous) pour l'action du groupe
  modulaire de \tei\ ${\rm Mod}(S)$.  Voir aussi \cite{Wol1,EL} pour
  une description de cette compactification par les applications
  harmoniques, et \cite{Wol2} pour la relation entre les applications
  harmoniques et les actions de groupes de $\pi_1S$ sur les arbres
  (réels). Kaimanovich et Mazur \cite{KM} ont montré que le bord de
  Poisson de nombreuses cha\^{\i}nes de Markov sur l'espace de \tei\
  est le bord de Thurston muni d'une mesure de probabilité adéquate.

\item[$\bullet$] {\bf La compactification naturelle de Morgan-Shalen}
  (voir \cite{MS1}). Soit $X$ une composante irrréductible de
  l'ensemble algébrique affine complexe $${\rm
    Hom}(\pi_1S,SL_2(\CC))//SL_2(\CC)\;,$$ qui est le quotient
  algébro-géométrique (voir \cite{MF}, et \cite{CS} pour une
  description simple) de l'ensemble algébrique affine complexe ${\rm
    Hom}(\pi_1S,SL_2(\CC))$ des morphismes de groupes de $\pi_1S$ dans
  $SL_2(\CC)$, par l'action par conjugaison au but de $SL_2(\CC)$.
  Pour tout $\ga$ dans $\pi_1S-\{1\}$, l'application polynomiale de
  ${\rm Hom}(\pi_1S,SL_2(\CC))$ dans $\CC$, définie par $\rho\mapsto
  {\rm trace}\;\rho(\ga)$, induit une application $f_\ga:X\ra \CC$,
  indépendante de la classe de $\ga$ dans l'ensemble $\C$ des classes
  de conjugaison d'éléments non triviaux de $\pi_1S$.  Si $\widehat X$
  est le compactifié d'Alexandrov de $X$, alors Morgan et Shalen
  \cite{MS1} montrent que l'application de $X$ dans $\widehat
  X\times\PP({\RR_+}^\C)$, définie par $x\mapsto
  (x,\RR_+^*(\log(|f_\ga(x)|+2))_{\ga\in\C})$, est un homéomorphisme
  sur son image, qui est ouverte dans son adhérence, et d'adhérence
  compacte.  L'espace de \tei\ de $S$ est homéomorphe, de manière
  équivariante par ${\rm Mod}(S)$, à une composante connexe de
  l'ensemble des points réels d'une composante irréductible $X$ de
  ${\rm Hom}(\pi_1S,SL_2(\CC))//SL_2(\CC)$. Morgan et Shalen
  obtiennent ainsi une compactification naturelle de $\ts$. Il n'est
  pas difficile de montrer avec cette définition qu'elle est isomorphe
  à la compactification de Thurston.  Morgan et Shalen \cite{MS1}
  montrent aussi (en version plus élaborée) que tout point du bord est
  de la forme $(\infty, \RR_+^*(\max\{0,-v(f_\ga)\})_{[\ga]\in\C})$ où
  $\infty$ est le point à l'infini de $\widehat X$ et $v:F^\times\ra
  \RR$ est une valuation réelle (pas forcément dis\-crète) sur le
  corps des fonctions $F$ de $X$, ainsi que de la forme $$(\infty,
  \RR_+^*(\min_{x\in T} d(x,\ga x) )_{[\ga]\in\C})$$ où $T$ est un
  arbre (réel) muni d'une action isométrique de $\pi_1S$, construit à
  partir de l'arbre de Bruhat-Tits de ${\rm SL}_2$ sur le corps valué
  $(F,v)$.  Le passage des feuilletages transversalement mesurés sur
  $S$ aux actions de $\pi_1S$ sur des arbres (réels) se fait par la
  notion d'arbre dual, développée dans \cite{MO,MS2}, voir aussi
  \cite{Sko,Ota,LP} par exemple.

\item[$\bullet$] {\bf La compactification naturelle de Bestvina}
  \cite{Bes} {\bf et de l'auteur} \cite{Pau1,Pau2}, qui est isomorphe
  à celle de Thurston. La présentation ci-dessous, utilisant la
  topologie de Gromov équivariante, correspond à une partie non
  publiée de la thèse de l'auteur \cite{Pau1}, mais il est possible de
  la retrouver aussi en partie dans la littérature actuelle (voir par
  exemple \cite{Ota,Kap}).

\item[$\bullet$] {\bf La compactification naturelle de Brumfiel}
  \cite{Bru}, voir aussi \cite{Bouz}. Elle utilise comme celle de
  Morgan-Shalen la description algébrique de l'espace de \tei, avec de
  plus la théorie du spectre réel des anneaux commutatifs. La
  compactification de Thurston est l'image par une surjection
  continue, équivariante par ${\rm Mod}(S)$, de cette
  compactification.
\item[$\bullet$] {\bf La compactification naturelle de Bonahon} par
  les courants géodésiques \cite{Bon1}, qui est isomorphe à celle de
  Thurston. Le plongement fonctionnel de Thurston de l'espace de
  Teichmüller, décrit ci-dessous, n'est sans doute pas le plus
  intéressant du point du vue analytique. Le plongement de Bonahon
  \cite{Bon1} semble plus approprié au calcul infinitésimal (voir
  \cite{Bon2,Bon3}). Ce plongement est défini de la manière suivante.
  Soit $\partial \wt S$ le bord à l'infini de $\wt S$ pour la métrique
  $\wt\sigma_0$, relevée par un rev\^etement universel $\wt S\ra S$
  d'une métrique hyperbolique fixée $\sigma_0$ sur $S$, muni de son
  action du groupe de revêtement $\pi_1 S$.  Remarquons qu'à
  homéomorphisme équivariant près, cet espace ne dépend pas de la
  métrique choisie $\sigma_0$.  Notons $\partial_2 \wt S$ l'espace
  topologique (métrisable séparable) localement compact des paires de
  points distincts de $\partial \wt S$, muni de l'action diagonale de
  $\pi_1 S$, et
  $$\M(\partial_2 \wt S)_{\pi_1S}$$ %
  l'espace vectoriel topologique (pour la topologie vague) des mesures
  de Radon $\pi_1S$-in\-va\-rian\-tes sur $\partial_2 \wt S$.  Pour
  tout élément $\sigma$ de $\ts$, le fibré unitaire tangent de la
  métrique relevée (d'un représentant) de $\sigma$ sur $\wt S$,
  quotienté par le changement de signe $v\mapsto -v$, s'envoie dans
  $\partial_2 \wt S$ par la fibration principale, de fibre $\RR$ pour
  l'action du flot géodésique, qui à un vecteur unitaire tangent
  associe la paire d'extrémités de la géodésique qu'il dirige.  Le
  relevé à $T^1\wt S$ de la mesure de Liouville de $(S,\sigma)$, qui
  est invariant par le flot géodésique, donne donc un élément
  $\lambda_\sigma$ de $\M(\partial_2 \wt S)_{\pi_1S}$.  Le plongement
  de Bonahon est l'application $\sigma\mapsto \lambda_\sigma$ de $\ts$
  dans $\M(\partial_2 \wt S)_{\pi_1S}$. La compactification de Bonahon
  est alors obtenue en projectifiant et en prenant l'adhérence.
\item[$\bullet$] Pour des surfaces particulières, il existe d'autres
  compactifications bien étudiées, telle la {\bf compactification de
    Maskit} de $\ts$ si $S$ est un tore privé d'un point (voir par
  exemple \cite{KS}, ses dessins et ses références).
\end{itemize}

\medskip Il n'est pas connu de l'auteur si deux compactifications de
\tei\ de $\ts$ ni si deux compactifications de Bers de $\ts$ (pour
deux choix de structures de surface de Riemann non isomorphes sur $S$)
sont isomorphes ou non. Il est montré dans \cite{Ker} (voir aussi
\cite{Mas}) que la compactification de Thurston n'est isomorphe à
aucune compactification de \tei. Certes, le bord de Teichm\"uller
s'envoie contin\^ument dans le bord de Thurston par l'application qui
à une forme différentielle quadratique holomorphe associe son
feuilletage transversalement mesuré horizontal (voir par exemple
\cite{Gar}). Mais cette application ne s'étend pas contin\^ument en un
isomorphisme de la compactification de \tei\ considérée sur celle de
Thurston. Il est montré dans \cite{KT} (voir aussi \cite{Bro}) que la
compactification de Thurston n'est pas isomorphe à au moins l'une des
compactifications de Bers (parce que pour au moins un des plongements
de Bers de l'espace de \tei, l'action du groupe modulaire ne s'étend
pas au bord).

\medskip
\begin{center}
+ --------------------------- +
\end{center}

\medskip Après cette disgression, entrons dans le vif du sujet.  Soit
$X$ un espace topologique localement compact (respectivement un espace
topologique localement compact muni d'une action par homéomorphismes
d'un groupe $G$). Une {\it compactification} (respectivement {\it
  compactification naturelle}) de $X$ est un couple $(i,Y)$ (ou par
abus $Y$) où $Y$ est un espace topologique compact, et $i:X\ra Y$ un
homéomorphisme sur son image, (par lequel on identifie $X$ et son
image dans $Y$) d'image ouverte et dense, (respectivement un tel
couple $(i,Y)$ tel que l'action de $G$ sur $X$ s'étende continuement
en une action par homéomorphismes de $G$ sur $Y$).  Deux
compactifications (respectivement compactifications naturelles) $Y,Y'$
de $X$ sont {\it isomorphes} si l'identité de $X$ dans $X$ s'étend
continuement en un homéomorphisme (respectivement un homéomorphisme
équivariant pour les actions de $G$) de $Y$ dans $Y'$.

Par exemple, les compactifications d'Alexandrov et de Stone-Cech (voir
par exemple \cite{Dug}) d'un espace localement compact $X$ sont
naturelles pour l'action du groupe $G$ de tous les homéomorphismes de
$X$.  Outre sa naturalité, l'intér\^et de la compactification de Thurston
des espaces de \tei\ est qu'elle apporte des informations
intéressantes sur les dégénérescences géométriques des structures
hyperboliques sur les surfaces.

\bigskip Soit $\HH^2_\RR$ le modèle du demi-espace supérieur du plan
hyperbolique réel, et ${\rm Isom_+}(\HH^2_\RR)$ le groupe des
isométries préservant l'orientation de $\HH^2_\RR$. L'action par
homographies de ${\rm SL}_2(\RR)$ sur $\HH^2_\RR$ induit un
isomorphisme (de groupes de Lie) de ${\rm PSL}_2(\RR)={\rm
  SL}_2(\RR)/\{\pm{\rm Id}\}$ avec ${\rm Isom_+}(\HH^2_\RR)$, par
lequel nous identifions ${\rm PSL}_2(\RR)$ et ${\rm
  Isom_+}(\HH^2_\RR)$.

Soit $S$ une surface compacte (sans bord) connexe orientée de classe
$\cinf$ privée d'un nombre fini de points, et de caractéristique
d'Euler strictement négative. En particulier, $S$ admet au moins une
{\it métrique hyperbolique} (i.e.~métrique riemannienne à courbure
constante $-1$) complète de volume fini (voir par exemple
\cite{FLP,Bus}).  Notons $\widetilde {S}\ra S$ un rev\^etement
universel de $S$, et $\pi_1S$ son groupe des automorphismes de
rev\^etement, qui, pour un choix fixé de point base dans $\wt S$,
s'identifie au groupe fondamental de $S$ en l'image $x_0$ de ce point
base dans $S$. Un élément de $\pi_1S-\{1\}$ sera dit {\it parabolique}
s'il est repr\'esentable par un lacet librement homotope dans tout
voisinage de l'un des bouts de $S$.

Notons ${\rm Diff}(S)$ le groupe des difféomorphismes $f$ de classe
$\cinf$ de $S$, muni de la topologie $\cinf$ (voir par exemple
\cite{Hir}), et ${\rm Diff}_0(S)$ le sous-groupe distingué ouvert des
difféomorphismes de $S$ isotopes à l'identité. Tout difféomorphisme
$f$ de $S$ induit un isomorphisme de groupes de $\pi_1(S,x_0)$ dans
$\pi_1(S,f(x_0))$, donc (par changement de point base, possible par
connexité) un isomorphisme $f_*:\pi_1S\ra\pi_1S$ bien défini modulo
conjugaison. De plus, si $f$ est isotope à l'identité, alors
$f_*:\pi_1S\ra\pi_1S$ est une conjugaison.

Rappelons qu'une métrique riemannienne sur $S$ est une section $\cinf$
du fibré, noté $\otimes^2T^*S$, des formes bilinéaires sur les espaces
tangents aux points de $S$. Notons $\wts$ l'ensemble des métriques
hyperboliques $\sigma$ sur $S$, complètes et de volume fini, muni de
la topologie induite par la topologie $\cinf$ sur
$\cinf(S,\otimes^2T^*S)$. Le groupe ${\rm Diff}(S)$ agit sur $\wts$
par l'action naturelle sur les champs de tenseurs sur $S$
$$(f,\sigma)\mapsto \left(f_*\sigma: x\in S\mapsto \{(u,v)\in
  (T_xS)^2\mapsto
  \sigma_{f^{-1}(x)}((T_xf)^{-1}(u),(T_xf)^{-1}(v))\}\right)\;.$$

L'{\it espace de Teichmüller} $\ts$ de $S$ est l'espace topologique
quotient de $\wts$ par l'action de ${\rm Diff}_0(S)$. Le {\it groupe
  modulaire de Teichmüller} (ou \og mapping class group \fg) de $S$
est le groupe discret quotient $${\rm Mod}(S)={\rm Diff}(S)/{\rm
  Diff}_0(S)\;.$$ %
L'action de ${\rm Diff}(S)$ sur $\wts$ induit une action du groupe
modulaire de Teichmüller de $S$ sur l'espace de Teichmüller de $S$.
Par abus, nous noterons de la même manière un élément $\sigma$ de
$\wts$ et son image dans $\ts$, ainsi qu'un élément $\sigma$ de $\ts$
et l'un de ses relevés dans $\wts$. Notons ${\rm Diff}_+(S)$ le
sous-groupe de ${\rm Diff}(S)$ préservant l'orientation et ${\rm
  Mod}_+(S)={\rm Diff}_+(S)/{\rm Diff}_0(S)$.

Par le théorème de Cartan-Hadamard (voir par exemple \cite{GHL}), pour
tout élément $\sigma$ de $\wts$, il existe une isométrie $i$
préservant l'orientation entre le relevé $\widetilde{\sigma}$ au
revêtement universel orienté de $S$ de la métrique $\sigma$ et le plan
hyperbolique réel $\HH^2_\RR$.  La conjugaison par une telle isométrie
$i$ fournit une {\it représentation fidèle et discrète} (i.e.~un
morphisme de groupes injectif d'image discrète)
$\rho=\rho_\sigma:\pi_1 S\ra {\rm Isom_+}(\HH^2_\RR)={\rm
  PSL}_2(\RR)$, appelée une {\it représentation d'holonomie} de
$\sigma$, et cette isométrie $i$ induit par passage au quotient une
isométrie préservant l'orientation de $(S,\sigma)$ sur la variété
hyperbolique quotient $\rho(\pi_1 S)\backslash \HH^2_\RR$.  La
représentation d'holonomie $\rho$ de $\sigma$ est indépendante du
choix de $i$ modulo conjugaison au but par un élément de ${\rm
  PSL}_2(\RR)$.  De plus, $\rho$ {\it préserve les paraboliques},
i.e.~envoie un élément parabolique de $\pi_1 S$  sur un élément
parabolique de ${\rm Isom_+}(\HH^2_\RR)$, ceci par la structure des
bouts des surfaces hyperboliques de volume fini (voir par exemple
\cite{Bus,Rat}).  Pour tout $f$ dans ${\rm Diff}(S)$, les
représentations $\rho_{f_*\sigma}$ et $\rho_\sigma\circ f_*$ sont
conjuguées. En particulier, si $\sigma$ et $\sigma '$ sont dans la
même orbite sous ${\rm Diff}_0(S)$, alors leurs représentations
d'holonomie sont conjuguées dans ${\rm PSL}_2(\RR)$ (car un
difféomorphisme isotope à l'identité induit une conjugaison sur
$\pi_1S$).

Considérons l'espace 
$$\rfd$$ 
des orbites, pour l'action par conjugaison au but de ${\rm PSL}_2
(\RR)$, des représentations fidèles et discrètes, préservant les
paraboliques, de $\pi_1 S$ dans ${\rm PSL}_2(\RR)$. Munissons-le de la
topologie quotient de la topologie de la convergence simple.
L'application $(f,\rho)\mapsto \rho\circ f_*$, pour $f$ dans ${\rm
  Diff}(S)$ et $\rho:\pi_1S\ra{\rm PSL}_2(\RR)$, induit donc une
action par homéomorphismes du groupe ${\rm Mod}(S)$ sur l'espace
$\rfd$.

On a une application de $\ts$ à valeurs dans cet espace, qui à
$\sigma$ associe la classe de conjugaison d'une représentation
d'holonomie de $\sigma$.

\bprop \label{pro:carcrepholo} (Voir par exemple \cite{Gol2}) 
Cette application est un homéomorphisme sur l'une des deux composantes
connexes de $\rfd$. Elle est équivariante pour l'action de ${\rm
  Mod}_+(S)$. \cqfd 
\eprop

L'un des points de la preuve de cette proposition est le suivant.

\blemm 
Une limite de représentations fidèles et discrètes,
préservant les paraboliques, de $\pi_1 S$ dans ${\rm PSL}_2(\RR)$ est
encore fidèle et discrète, préservant les paraboliques.  
\elemm

\dem Soit $(\rho_i)_{i\in\NN}$ une suite de telles représentations,
qui converge vers $\rho$. Nous renvoyons par exemple à \cite{GM}
pour le caractère fidèle et la discrétude de $\rho$.  Comme un
élément $\pm A$ de ${\rm PSL}_2(\RR)$ est parabolique ou
l'identité si et seulement si $|{\rm tr}\;A|=2$, le fait que $\rho$
préserve les paraboliques en découle par continuité de la trace.  
\cqfd

\medskip Notons $\C$ l'ensemble des classes de conjugaison des
éléments non paraboliques de $\pi_1S-\{1\}$, qui s'identifie avec
l'ensemble des classes d'homotopie libre de lacets de $S$,
homotopiquement non triviaux, et non librement homotopables dans tout
voisinage d'un des bouts de $S$.
%, ou de manière
%équivalente, l'ensemble des classes de conjugaison des éléments $\ga$
%de $\pi_1S-\{1\}$ tels que si $\rho$ est une représentation
%d'holonomie d'un élément quelconque de $\ts$, alors $\rho(\ga)$ soit
%parabolique).  
Pour tout $\alpha$ dans $\C$ et tout $\sigma$ dans $\ts$, on note
$\ell_\sigma(\alpha)$ la longueur de l'unique géodésique fermée pour
$\sigma$ dans la classe d'homotopie libre $\alpha$ (voir par exemple
\cite{Bus} pour l'existence et l'unicité de cette géodésique).  Donc
$\ell_\sigma(\alpha)>0$ est égale à la {\it distance de translation}
$$\ell_\sigma(\alpha)=\inf_{x\in\HH^2_\RR} d(x,\rho(\ga)x)$$
de $\rho(\ga)$ dans $\HH^2_\RR$, pour $\rho$ une représentation
d'holonomie de $\sigma$ et pour n'importe quel $\ga\in\pi_1S$ dans la
classe $\alpha$. Il n'est pas difficile de montrer que, pour tout
$\ga$ dans $\C$, l'application $\sigma\mapsto\ell_\sigma(\alpha)$ de
$\ts$ dans $\RR$ est continue (un moyen rapide est d'utiliser, outre
la définition de la topologie de $\ts$, le lemme de fermeture d'Anosov
\cite{Ano}, enfin une version très simplifiée disant qu'une courbe
fermée presque géodésique est à distance de Hausdorff petite d'une
géodésique fermée, voir par exemple \cite{BH}).

En notant $\RR_+=[0,+\infty[$, on a donc une application
$$\begin{array}[b]{cccc}
\ell:&\ts & \ra & \RR_+\mbox{}^{\C}-\{0\}\\
&\sigma & \mapsto & \left(\ell_\sigma(\alpha)\right)_{\alpha\in\C}
\end{array}\;.
$$ %
Notons $\pi:(\RR_+\mbox{}^{\C}-\{0\})\ra \PP(\RR_+\mbox{}^{\C})=
\left(\RR_+\mbox{}^{\C}-\{0\}\right)/\sim$ la projection canonique
pour la relation d'équivalence définie par $(x_\alpha)_ \alpha\sim
(y_\alpha)_\alpha$ si et seulement s'il existe $\lambda>0$ tel que
$x_\alpha=\lambda y_\alpha$ pour tout $\alpha$. On munit l'espace
produit $\RR_+\mbox{}^{\C}$ de la topologie produit, et l'espace
quotient $\PP(\RR_+\mbox{}^{\C})$ de la topologie quotient. Le groupe
des difféomorphismes de $S$ agit sur $\C$, par l'application qui à
$(f,\alpha)\in{\rm Diff}(S)\times\C$ associe la classe d'homotopie
libre de $f(\ga)$ pour n'importe quel lacet $\ga$ dans la classe
d'homotopie libre $\alpha$, et l'action de ${\rm Diff}_0(S)$ est
triviale. Par conséquent ${\rm Mod}(S)$ agit sur $\C$, donc agit par
homéomorphismes sur $\RR_+\mbox{}^{\C}$ par permutation des
coordonnées, et donc agit par homéomorphismes sur
$\PP(\RR_+\mbox{}^{\C})$ par passage au quotient.

\btheo \label{theo:benouaisquoi} (W.~Thurston \cite{Thu0})
L'application $\pi\circ\ell:\ts\ra\PP(\RR_+\mbox{}^{\C})$ est un
homéomorphisme sur son image, qui est relativement compacte, et
ouverte dans son adhérence. Elle est équivariante pour les actions de
${\rm Mod}(S)$.  
\etheo

Ainsi, le couple $(\pi\circ\ell,\,\ov{\pi\circ\ell(\ts)}\;)$ est une
compactification de $\ts$, appelée la {\it compactification de
  Thurston} de l'espace de Teichmüller, qui est naturelle pour
l'action de ${\rm Mod}(S)$.

Dans cet énoncé, nous pourrions remplacer l'ensemble $\C$ par
l'ensemble des classes d'homotopie libre de lacets simples, non
homotopes à $0$, ni homotopables dans tout voisinage d'un point enlevé
de $S$, comme dans \cite{FLP}, ou par l'ensemble de toutes les classes
de conjugaison non triviales de $\pi_1 S$, i.e.~de toutes les classes
d'homotopie libre de lacets non homotopes à $0$. Les résultats de
\cite{FLP} impliquent les deux autres versions dans ce qui suit.

\bigskip 
\noindent {\bf Preuve du théorème \ref{theo:benouaisquoi}. } 
Nous n'utiliserons pas l'approche originale de Thurston \cite{Thu0}
exposée complètement dans \cite{FLP}, mais celle de
\cite{Bes,Pau1,Pau2}, aussi présentée sous une forme un peu différente
dans \cite{Ota}. Voir \cite{Par} pour de profondes généralisations.

Nous admettrons (voir par exemple \cite{FLP,Bus}) ce qui concerne la
description interne de l'espace de Teichmüller, résumé dans le fait
suivant.
%,qui a essentiellement fait l'objet des exposés de L.~Funar et
%V.~Sergiescu.

\btheo \label{theo:admisflp} %
Il existe un nombre fini $\ga_1,\dots,\ga_N$ d'éléments de $\C$ (que
l'on peut même choisir représentables par des lacets simples) tels que
l'application $\ts\ra{\RR_+}^N$ définie par
$$\sigma\mapsto (\ell_\sigma(\ga_1),\dots,\ell_\sigma(\ga_N))$$ soit
continue, injective, propre, donc un homéomorphisme sur son image. En
particulier, $\ts$ est métrisable, séparable, et localement compact.
De plus $\ts$ est homéomorphe à $\RR^{6g-6+2b}$ si $g$ est le genre de
$S$ et $b$ le nombre de points enlevés de $S$.  
\cqfd 
\etheo

En particulier, l'application continue $\ell:\ts\ra
\RR_+\mbox{}^{\C}-\{0\}$ est un homéomorphisme propre sur son image.

\bprop 
L'application $\pi\circ\ell:\ts\ra\PP(\RR_+\mbox{}^{\C})$ est
un homéomorphisme sur son image, équivariant pour les actions de
${\rm Mod}(S)$. 
\eprop

\medskip \dem (Voir \cite[exposé 7]{FLP}.) Pour montrer l'injectivité
de $\pi\circ\ell$, supposons par l'absurde qu'il existe des éléments
$\sigma$ et $\sigma'$ dans $\ts$ et $t>1$ tels que $\ell_{\sigma'} =
t\ell_\sigma$. Soient $a,b$ dans $\pi_1 S$ tels que
$\rho_\sigma(a),\rho_\sigma(b),\rho_\sigma(ab), \rho_\sigma(a^{-1}b)$
soient des éléments hyperboliques, avec
$\rho_\sigma(a),\rho_\sigma(b)$ d'axes de translation
(transversalement) sécants, et tels que
$\rho_\sigma(a),\rho_\sigma(b)$ se relèvent en des matrices $A,B$ de
${\rm SL}_2(\RR)$ telles que les traces ${\rm tr}(A),{\rm tr}(B),{\rm
  tr}(AB),{\rm tr}(A^{-1}B)$ soient strictement positives, et de
m\^eme en remplaçant $\sigma$ par $\sigma'$.

Ceci est possible pour $\sigma$, par exemple en fixant $a$ d'image par
$\rho_\sigma$ hyperbolique, et en prenant pour $b$ une puissance $N$
suffisamment grande d'un élément de $\pi_1S$ d'image par $\rho_\sigma$
hyperbolique, dont le point attractif de l'image par $\rho_\sigma$
soit dans un voisinage $U_+$ suffisamment petit de celui de $a$, et le
point fixe répulsif soit dans un voisinage $U_-$ suffisamment petit,
et du bon côté, de celui de $a$. Ceci est possible par densité des
couples de points fixes d'éléments hyperboliques de
$\rho_\sigma(\pi_1S)$ dans l'ensemble des couples de points du bord à
l'infini de $\HH^2_\RR$ (voir par exemple \cite{Gro} dans un cadre
bien plus général). Rappelons que les actions de $\rho_\sigma(\pi_1S)$
et $\rho_{\sigma'}(\pi_1S)$ sur le cercle à l'infini $\SS^1_\infty$ de
$\HH^2_\RR$ sont conjuguées, par l'extension continue à l'infini du
relevé à $\HH^2_\RR$ de l'application identité de $S$ dans $S$ par les
rev\^etements universels riemanniens $\HH^2_\RR\ra (S,\sigma)$ et
$\HH^2_\RR\ra (S,\sigma')$ (voir par exemple \cite{ET} pour une preuve
de cette extension y compris en dimension supérieure, souvent
faussement attribuée à Mostow qui s'en sert dans \cite{Mos}, et qui
dans le cas de $\HH^2_\RR$ devait \^etre déjà connue). Donc en
choisissant $N$ suffisamment grand et $U_+,U_-$ suffisamment petits, on
obtient bien les m\^emes propriétés pour $\sigma'$.

\blemm 
La distance de translation $\ell(A)$ dans $\HH^2_\RR$ d'un
élément hyperbolique $A$ de ${\rm SL}_2(\RR)$ vérifie
$$\cosh \frac{\ell(A)}{2} =  \left|\frac{{\rm tr}(A)}{2}\right|\;.$$
Si $A$ et $B$ sont deux éléments de ${\rm SL}_2(\RR)$, alors
$$ {\rm tr}(AB)+{\rm tr}(A^{-1}B)={\rm tr}(A){\rm tr}(B)\;.$$
\elemm

\dem Le premier résultat découle du fait que $A$ est conjuguée à une
matrice de la forme $\pm \left(\begin{array}{cc} \lambda & 0 \\ 0&
    \lambda^{-1}\end{array}\right)$ où $\lambda>1$, et que la distance
de translation de cette matrice vaut $2\log \lambda$ (par unicité
d'une géodésique translatée, son axe de translation est l'axe
vertical).

Pour montrer le second résultat, le théorème de Cayley-Hamilton
implique que $$A^2 - {\rm tr}(A) A + {\rm Id}=0\;,$$ donc $A+A^{-1} =
{\rm tr}(A){\rm Id}$, d'où le résultat en multipliant par $B$ à droite
et en prenant la trace.  
\cqfd

\medskip %
Maintenant, comme $\ell_{\sigma'} = t\ell_\sigma$, en notant $\ell_1,
\ell_2,\ell_3,\ell_4$ la distance de translation de
$\rho_\sigma(a),\rho_\sigma(b),\rho_\sigma(ab)$,
$\rho_\sigma(a^{-1}b)$ respectivement, on déduit du lemme précédent
que
$$\cosh (\frac{\ell_3}{2}) +
\cosh (\frac{\ell_4}{2})=2\cosh (\frac{\ell_1}{2})\cosh
(\frac{\ell_2}{2})=\cosh (\frac{\ell_1+\ell_2}{2})+\cosh
(\frac{\ell_1-\ell_2}{2})\;.$$

Le fait suivant est laissé en exercice au lecteur.

\blemm %
Soient $u,v,w,x,t$ des nombres réels strictement positifs, 
avec $t>1$, tels que
$$\cosh u +\cosh v = \cosh w +\cosh x, \;\;\;
\cosh tu +\cosh tv = \cosh tw +\cosh tx\;.$$ 
Alors $\{u,v\}=\{w,x\}$.
\cqfd
\elemm

Maintenant, on aurait, quitte à échanger $ab$ et $a^{-1}b$, l'égalité
$\ell_3=\ell_1+\ell_2$. Comme les axes de translation (des images par
$\rho_\sigma$) de $a$ et de $b$ se coupent (transversalement), ceci
contredit la minimalité de la distance de translation $\ell_3$.
L'injectivité de $\pi\circ\ell$ en découle.

\medskip Supposons qu'il existe une suite $(\sigma_i)_{i\in\NN}$ dans
$\ts$, un élément $\sigma$ dans $\ts$ et une suite de réels
strictement positifs $(t_i)_{i\in\NN}$, tels que $t_i\ell_{\sigma_i}$
converge vers $\ell_\sigma$ dans $\RR_+\mbox{}^{\C}$.  Si
$(t_i)_{i\in\NN}$ convergeait vers $+\infty$, alors par la propreté
dans le théorème \ref{theo:admisflp}, quitte à extraire, la suite
$(\sigma_i)_{i\in\NN}$ convergerait vers un élément $\sigma'$ de $\ts$
tel que $\ell_{\sigma'}=0$, ce qui est impossible.  Donc quitte à
extraire, $t_i$ converge vers $t\geq 0$. Soit $s$ la borne inférieure
des $\ell_{\sigma}(\alpha)$ pour $\alpha$ dans $\C$, qui est
strictement positive, par le lemme de Zassenhauss, voir par exemple
\cite{Rag}, (le lemme de Zassenhauss n'est qu'un cas très particulier
(et antérieur) pour les espaces localement symétriques du lemme de
Margulis, voir par exemple \cite{BK}). Comme le volume de $\sigma_i$
est constant (égal à $2\pi|\chi(S)|$), la borne inférieure des
$\ell_{\sigma_i}(\alpha)$ pour $\alpha$ dans $\C$ est majorée
uniformément en $i$, donc $t>0$.  Par conséquent, la suite
$(\ell_{\sigma_i})_{i\in\NN}$ converge vers $\frac{1}{t}\ell_\sigma$.
Par la propreté dans le théorème \ref{theo:admisflp}, quitte à
extraire, $(\sigma_i)_{i\in\NN}$ converge vers un élément $\sigma'$
dans $\ts$. Par continuité, nous avons
$\ell_{\sigma'}=\frac{1}{t}\ell_\sigma$.  Par l'injectivité montrée
ci-dessus, nous avons $t=1$, et $(\sigma_i)_{i\in\NN}$ converge vers
$\sigma'=\sigma$.

Ceci implique que l'application continue injective $\pi\circ\ell$ est
un homéomorphisme sur son image.  La naturalité est évidente.
\cqfd

\bigskip %
Soit $\Ga$ un groupe, et $\E$ un ensemble d'espaces métriques munis
d'une action isométrique de $\Ga$. Pour tout $X$ dans $\E$, pour toute
partie finie $K$ de $X$, pour tout $\epsilon >0$ et toute partie finie
$P$ de $\Ga$, notons $\V_{\epsilon,P,K}(X)$ l'ensemble des éléments
$X'$ de $\E$ tels qu'il existe une partie finie $K'$ de $X'$, et une
relation $\R\subset K\times K'$, dont les deux projections sur $K$ et
$K'$ sont surjectives, telles que
$$\forall\; x,y\in K,\;\forall\; x',y'\in K',\;\forall\; \ga\in
P,\;\;{\rm si}\;\;x\,\R\, x'\;\,{\rm et}\;\,y\,\R\, y'\;\;
{\rm alors}\;\; |d(x,\ga y)-d(x',\ga y')|<\epsilon\;.$$

%\begin{enumerate}
%\item[$\bullet$] $\forall\; x\in K,\forall\; x'\in K',\forall\; \ga\in
%  P,\;\;{\rm si}\;\;x\,\R\, x'\;\;{\rm et}\;\;\ga x\in K\;\;\;{\rm
%    alors}\;\;\; \ga x'\in K'\;\;{\rm et}\;\;\ga x\,\R\,\ga x'$,
%\item[$\bullet$] 
%$\forall\; x,y\in K,\;\forall\; x',y'\in K',\;\forall\; \ga\in
%  P,\;\;{\rm
%    si}\;\;x\,\R\, x'\;,\;y\,\R\, y'\;\;{\rm alors}\;\;
%|d(\ga x,y)-d(\ga x',y')|<\epsilon$.
%\end{enumerate}

\noindent %
La {\it topologie de Gromov équivariante} sur $\E$ est la
topologie sur $\E$ (qui existe, voir \cite{Pau2}) dont l'ensemble des
$\V_{\epsilon,P,K}(X)$ est un système fondamental de voisinages de
$X$, pour tout $X$ dans $\E$.

Par exemple, l'ensemble $\wts$ s'identifie, par l'application de
relèvement $\sigma\mapsto\wt\sigma$, à l'ensemble des métriques
hyperboliques complètes sur le revêtement universel $\wt S$ de $S$,
invariantes par l'action du groupe de revêtement $\pi_1S$, telles que
tout élément parabolique de $\pi_1S$ agisse par une isométrie
parabolique sur $\wt S$.

\bprop \label{pro:coincidtopos}
(F.~Paulin \cite[Prop.~6.2]{Pau2}) La topologie quotient sur $\ts$ de 
la topologie de Gromov sur $\wts$ co\"incide avec la topologie
usuelle de $\ts$.
\eprop

\dem Il est immédiat que la topologie de Gromov équivariante quotient
est moins fine que la topologie usuelle. En effet, supposons que deux
métriques hyperboliques complètes de volume fini sur $S$ soient
proches. Alors elles ont des voisinages de bouts proches qui sont
isométriques, par la structure des bouts des surfaces hyperboliques.
Et sur le compact restant, les distances riemanniennes sont proches.
Par passage au revêtement universel, on montre bien que pour toute
partie finie $K$ de $\wt S$, pour toute partie finie $P$ de $\pi_1S$
et tout $\epsilon>0$, si $\sigma'$ est suffisamment proche de $\sigma$
dans $\ts$, alors $(\wt S,\wt\sigma')$ appartient à
$\V_{\epsilon,P,K}((\wt S,\wt\sigma))$.

Réciproquement, montrons que pour tout $\epsilon >0$ et tout $\alpha$
dans $\C$, si $(\wt S,\wt\sigma')$ est suffisamment proche de $(\wt
S,\wt\sigma)$ pour la topologie de Gromov équivariante, alors
$|\ell_{\sigma'}(\alpha)-\ell_\sigma(\alpha)|<\epsilon$. Comme
$\ell:\ts\ra \RR_+\mbox{}^{\C}-\{0\}$ est un homéomorphisme sur son
image, ceci montrera le résultat.  En effet, soit $\ga\in\pi_1 S$ un
représentant de $\alpha$, et $x$ un point de l'axe de translation de
$\rho_\sigma(\ga)$. Posons $K=\{x\}$ et $P=\{\ga,\ga^2\}$. Si $(\wt
S,\wt\sigma')$ appartient à $\V_{\eta,P,K}((\wt S,\wt\sigma))$, alors
il existe un point $x'$ de $\wt S$ tel que, en notant $d_{\wt\sigma'}$
la distance de la variété riemanienne $(\wt S,\wt\sigma')$, on ait
$$|d_{\wt\sigma'}(x',\ga x')-\ell_\sigma(\alpha)|<\eta$$
et 
$$|d_{\wt\sigma'}(x',\ga x')+d_{\wt\sigma'}(\ga x',\ga^2
x')-d_{\wt\sigma'}( x',\ga^2 x')|<3\eta\;.$$ %
Par de petites estimées de géométrie hyperbolique (en fait valable
dans n'importe quel espace métrique CAT$(-1)$, comme les arbres
(réels) définis ci-dessous, voir par exemple \cite{BH}), un chemin
géodésique par morceaux, dont la longueur totale est proche (par
rapport à la longueur des morceaux) de la distance entre ses
extrémités, est proche de la géodésique entre ses extrémités. Ceci
implique en itérant par $\gamma$ que le chemin géodésique par morceaux
$c=\bigcup_{n\in\ZZ}\;[\ga^n x',\ga^{n+1}x']$ a ses deux points à
l'infini distincts, et est proche de la géodési\-que $\overline c$
entre ces deux points. Comme $c$ est invariante par $\ga$, il en est
de m\^eme de $\overline c$. Ceci implique (ainsi que dans n'importe
quel espace métrique CAT$(-1)$, comme les arbres (réels)), lorsque
$\eta$ est assez petit, que $\rho_{\wt\sigma'}(\ga)$ est hyperbolique,
que le point $x'$ est proche de l'axe de translation de
$\rho_{\wt\sigma'}(\ga)$, et donc que la distance de translation de
$\rho_{\wt\sigma'}(\ga)$ est proche de celle de
$\rho_{\wt\sigma}(\ga)$.  \cqfd

\bigskip Le résultat clef dans l'étude des dégénérescences de
structures hyperboliques sur les surfaces est le théorème
\ref{theo:dichotomiequoi} suivant, voir \cite{Bes,Pau2}. La preuve de
ce théorème est essentiellement tirée de \cite{Pau4} (qui donne un
résultat bien plus général pour des actions proprement discontinues de
groupes de type fini sur des espaces symétriques de type non
compact). Donnons d'abord quelques définitions.

\medskip
Un {\it arbre (réel)} est un espace métrique $T$ {\it géodésique}
(i.e.~pour tous $x,y$ dans $T$, il existe une isométrie (appelée,
ainsi que son image, un {\it segment géodésique} entre $x$ et $y$)
$[a,b]\ra T$ envoyant $a$ sur $x$ et $b$ sur $y$), et {\it
  $0$-hyperbolique} (i.e.~pour tous $x,y,z$ dans $T$, pour tous les
segment géodésiques $[x,y]$, $[y,z]$, $[x,z]$ entre respectivement $x$
et $y$, $y$ et $z$, $x$ et $z$, le segment $[x,y]$ est contenu dans la
réunion $[y,z]\cup [x,z]$).

\medskip Rappelons (voir par exemple \cite{Tit}) que dans un arbre
(réel) $T$, toute isométrie $g$ ou bien admet un point fixe (et on dit
que $g$ est {\it elliptique} dans $T$) ou bien $T$ possède une unique
droite géodésique (appelée l'{\it axe de translation} de $g$) sur
laquelle $g$ agit par translation de distance
$\ell_T(g)=\min\{d(x,gx)\;:x\in T\}>0$ (appelée la {\it distance de
  translation} de $g$ dans $T$). Cet axe de translation est construit
ainsi : pour tout $x$ dans $T$, les segments $[x,gx]$ et $[gx,g^2x]$
se rencontrent en un segment $[gx,gu]$ pour $u\in[x,gx]$ par
$0$-hyperbolicité, et l'axe de translation est alors la réunion des
segments d'intérieurs disjoints $[g^ku, g^{k+1}u]$. Toute isométrie de
$T$ conjuguée à $g$ possède la même distance de translation que $g$.

Une action de $\pi_1 S$ sur un arbre (réel) $T$ est dite {\it à petits
  stabilisateurs d'arête} si c'est une action isométrique, sans point
fixe global, dont les fixateurs des segments géodésiques non réduits à
un point sont triviaux ou infinis cycliques, et telle que tout élément
parabolique de $\pi_1 S$ agisse par une isométrie elliptique sur
$T$. Ceci est une définition adaptée à la classe des groupes
(munis d'une classe de sous-groupes) que nous considérons.

Voir par exemple \cite{BT,Tit,MS1,MS2,Sha1,Sha2,Mor,Bes2,Pau5,Chi,FJ}
pour plus d'informations sur les arbres (réels).

\medskip
Pour $X$ un espace métrique de distance $d$, et $\epsilon >0$, 
notons $\epsilon X$ l'espace métrique $(X,\epsilon d)$.

\btheo \label{theo:dichotomiequoi}  %
Soit $(\sigma_i)_{i\in\NN}$ une suite dans $\ts$. Alors, quitte à
extraire, ou bien $(\sigma_i)_{i\in\NN}$ converge vers un élément
$\sigma$ de $\ts$, ou bien il existe une suite
$(\epsilon_i)_{i\in\NN}$ de nombres réels strictement positifs tendant
vers $0$, et un arbre (réel) $T$ muni d'une action à petits
stabilisateurs d'arête, de $\pi_1S$, tels que la suite
$(\epsilon_i\HH^2_\RR, \rho_{\sigma_i})$ converge vers $T$ pour la
topologie de Gromov équivariante.  
\etheo

\dem
Soit $\{s_1,\dots,s_p\}$ une partie génératrice finie de
$\pi_1S$. 

\blemm \label{lem:bestvina}(Voir \cite{Bes}) %
Soit $\rho:\pi_1S\ra{\rm PSL}_2(\RR)$ une représentation
fidèle et discrète. Alors il existe un point $x$ de $\HH^2_\RR$
qui minimize la distance $f_\rho(x) = \max_{1\leq j\leq p}
d(x,\rho(s_j)x)$.  
\elemm

\dem %
Soit $(x_i)_{i\in\NN}$ une suite dans $\HH^2_\RR$ telle que la suite
$f_\rho(x_i)$ converge vers la borne inférieure de $f_\rho$ sur
$\HH^2_\RR$.  Si $x_i$ reste dans un compact de $\HH^2_\RR$, alors,
quitte à extraire, $x_i$ converge vers un point $x$ qui convient.
Sinon, quitte à extraire, $x_i$ converge vers un point à l'infini $x$,
qui est fixe par $\rho(\pi_1S)$, ce qui est impossible, car le
stabilisateur dans ${\rm PSL}_2(\RR)$ d'un point à l'infini de
$\HH^2_\RR$ est résoluble, et $\pi_1S$ contient au moins un groupe
libre de rang $2$.  \cqfd

\medskip
Notons $\rho_{i}=\rho_{\sigma_i}$, et 
$$\lambda_i=\min_{x\in\HH^2_\RR}f_{\rho_{i}}(x)\;.$$ Fixons $*_i$ un
point de $\HH^2_\RR$ tel que $f_{\rho_i}(*_i)=\lambda_i$.  Alors,
quitte à extraire, ou bien la suite $(\lambda_i)_{i\in\NN}$ est
bornée, disons par $M$, ou bien elle converge vers $+\infty$.

Dans le premier cas, quitte à conjuguer $\rho_i$, nous pouvons
supposer que $*_i$ est un point fixé $*_0$ de $\HH^2_\RR$.  Par
compacité de l'ensemble des isométries de $\HH^2_\RR$ bougeant $*_0$
d'une distance au plus $M$, quitte à extraire, les $\rho_i(s_j)$ pour
$1\leq j\leq p$ convergent. Donc la suite $(\rho_i)_{i\in\NN}$
converge, et la proposition \ref{pro:carcrepholo} montre que la suite
$(\sigma_i)_{i\in\NN}$ converge (il suffit aussi de remarquer que la
suite $(\ell_{\sigma_i})_{i\in\NN}$, qui ne dépend que de
$(\rho_i)_{i\in\NN}$, converge, et d'utiliser le commentaire suivant
le théorème \ref{theo:admisflp}).

\bigskip
Supposons maintenant que la suite $(\lambda_i)_{i\in\NN}$ converge
vers $+\infty$.

Posons $(X_i,d_i)=\frac{1}{\lambda_i}\HH^2_\RR$. Soit $\omega$ un
ultrafiltre sur $\NN$, plus fin que le filtre de Fréchet des
complémentaires des parties finies (voir par exemple \cite{Bou}).
Pour toute suite $(x_i)_{i\in \NN}$ dans un espace topologique, notons
$\lim_\omega x_i$ une limite de cette suite pour le filtre image de
$\omega$ (voir par exemple \cite{Bou}), si elle existe (ce qui est le
cas si la suite $(x_i)_{i\in \NN}$ est relativement compacte; cette
limite est alors unique dans une partie compacte (donc séparée)
contenant la suite, et c'est une valeur d'adhérence de la suite).
Posons
$$X_\infty=\{x=(x_i)_{i\in\NN}\in\prod_{i\in\NN} X_i\;:\; \lim_\omega
d_i(x_i,*_i) <+\infty\}\;,$$
muni de la pseudo-distance
$$d_\infty(x,y)=\lim_\omega \;d_i(x_i,y_i)\;.$$ %
Notons $(X_\omega,d_\omega,*_\omega)$ l'espace métrique pointé
quotient de $(X_\infty,d_\infty,(*_i)_{i\in\NN})$ par la relation
d'équivalence définie par $x\sim y$ si $x$ et $y$ sont à
pseudo-distance nulle.  Pout tout $\gamma$ dans $\pi_1S$, l'isométrie
$\rho_i(\ga)$ de $\frac{1}{\lambda_i}\HH^2_\RR$ bouge le point $*_i$
d'une distance bornée (au plus $1$ si $\ga\in\{s_1,\dots,s_p\}$) pour
la distance $d_i$, donc au plus $n$ si $\ga$ s'écrit comme un mot
de longueur $n$ en $s_1,\dots, s_n$ et leurs inverses. Donc l'action
diagonale $(\ga,(x_i)_{i\in\NN})\mapsto (\rho_i(\ga)x_i)_{i\in\NN}$ de
$\pi_1S$ sur $X_\infty$ préserve $X_\infty$, est isométrique pour la
pseudo-distance $d_\infty$, donc induit une action isométrique
$\rho_\omega$ de $\pi_1S$ sur $X_\omega$.

Un petit calcul de géométrie hyperbolique montre que tout côté d'un
triangle géodésique de $\HH^2_\RR$ est à distance au plus
$\log(1+\sqrt{2})$ de la réunion des deux autres côtés. Donc tout côté
d'un triangle géodésique de $(X_i,d_i)$ est à distance au plus
$\log(1+\sqrt{2})/\lambda_i$ de la réunion des deux autres côtés.
C'est alors un exercice (voir par exemple \cite{KL}) de montrer que
l'espace métrique $(X_\omega,d_\omega)$ est géodésique, et
$0$-hyperbolique, donc est un arbre (réel). Le m\^eme argument montre
que toute limite, pour la topologie de Gromov équivariante, d'une
suite d'actions isométriques d'un groupe fixé sur des arbres (réels),
est encore une action isométrique de ce groupe sur un arbre (réel). Il
n'est pas difficile non plus de montrer que, pour la topologie de
Gromov équivariante,
$$\lim_\omega \;(X_i,d_i,\rho_i)= (X_\omega,d_\omega,\rho_\omega)\;.$$

\blemm\label{lem:limptstabar} %
L'action de $\pi_1S$ sur l'arbre (réel)
$X_\omega$ est à petits stabilisateurs d'arête.  
\elemm

\dem Par continuité, le point $*_\omega$ est un point bougé le moins
par les générateurs, et il est bougé au moins d'une distance de $1$
par l'un des générateurs, donc l'action n'a pas de point fixe global.

Supposons par l'absurde qu'il existe deux points distincts $x$ et $y$
de $X_\omega$, qui soient fixés par un sous-groupe non trivial, non
infini cyclique de $\pi_1 S$. Alors ils sont fixés par des éléments
$\alpha,\beta$ qui engendrent un sous-groupe libre sur $\alpha,\beta$
dans $\pi_1 S$.

Soit $\eta>0$, $P$ l'ensemble des mots en $\alpha,\beta$ de longueur
au plus $6$, et $K=\{x,y\}$. Pour tout $i$ suffisamment grand (pour
l'ultrafiltre $\omega$), le triplet $(X_i,d_i,\rho_i)$ appartient à
$\V_{\epsilon,P,K}((X_\omega,d_\omega,*_\omega))$, donc si $x$ est en
relation (voir la définition de la topologie de Gromov équivariante)
avec $x_i\in X_i$ et $y$ avec $y_i\in X_i$, alors pour tout $\ga$ dans
$P$,
$$d_{\HH^2_\RR}(x_i,\ga x_i)\leq \lambda_i \eta\;\;\;{\rm et}\;\;\;
d_{\HH^2_\RR}(x_i,y_i)\geq \lambda_i(d_\omega(x,y)- \eta)\;.$$ En
particulier, si $\eta$ a été choisi assez petit devant
$d_\omega(x,y)$, alors $\alpha$ et $\beta$ agissent presque comme des
translations le long d'un segment de longueur $6\lambda_i\eta$ centré
au milieu $m$ du segment géodésique $[x_i,y_i]$. Ceci implique que les
éléments $[\alpha,\beta], [\alpha,\beta^2]$, qui engendrent un groupe
libre, bougent le point $m$ d'une distance au plus $\eta$.  Si $\eta$
est assez petit, ceci contredit le lemme de Zassenhauss.

Supposons par l'absurde qu'il existe un élément parabolique $\gamma$
de $\pi_1S$ qui ne soit pas elliptique dans $X_\omega$. Soit $x$ un
point de l'axe de translation de $\gamma$ dans $X_\omega$. Soient
$\eta>0$, $K=\{x\}$ et $P=\{\ga,\ga^2\}$. Alors, pour $i$ suffisamment
grand (pour $\omega$), l'élément $(X_i,d_i,\rho_i)$ appartient à
$\V_{\eta,P,K}((X_\omega,d_\omega,\rho_\omega))$. Il existe donc $x_i$
dans $X_i$ tel que
$$|d_i(x_i,\rho_i(\ga) x_i)-\ell_{X_\omega}(\ga)|<\eta$$
et
$$|d_i(x_i,\rho_i(\ga) x_i)+d_i(\rho_i(\ga)x_i,\rho_i(\ga^2) x_i
)-d_i( x_i,\rho_i(\ga^2) x_i)|<3\eta\;.$$ %
Si $\eta$ est choisi petit devant $\ell_{X_\omega}(\ga)$, ceci
implique, comme dans la fin de la preuve de la proposition
\ref{pro:coincidtopos}, que $\rho_i(\ga)$ est hyperbolique, ce 
qui est une contradiction.  
\cqfd

\medskip
Ceci termine la preuve du théorème \ref{theo:dichotomiequoi}.
\cqfd

\medskip Si $T$ est un arbre (réel), muni d'une action isométrique du
groupe $\pi_1S$, telle que tout élément parabolique de $\pi_1S$ soit
elliptique dans $T$, alors notons $\ell_T\in \RR_+\mbox{}^{\C}$
l'application $\alpha\mapsto \ell_T(g)$ où $g\in\pi_1S$ appartient à
la classe $\alpha$.

Si $\lim_\omega (X_i,d_i,\rho_i)$ est un arbre (réel) $T$ muni d'une
action de isométrique de $\pi_1S$, alors le même argument que pour
montrer la réciproque de la proposition \ref{pro:coincidtopos} montre
que la suite $(\frac{1}{\lambda_i} \ell_{\sigma_i}(\alpha))
_{\alpha\in\C}$ converge vers l'élément $\ell_T$ dans ${\RR_+}^\C$.

Par métrisabilité de $\PP(\RR_+\mbox{}^{\C})$, ceci montre que l'image
de $\pi\circ\ell$ est d'adhérence compacte. Montrons qu'elle est
ouverte dans son adhérence.

\medskip
Il découle par exemple de \cite{CM} que si $(T_i)_{i\in\NN}$ est
une suite d'arbres (réels) munis d'une action isométrique de $\pi_1S$
à petits stabilisateurs d'ar\^ete, si la suite
$(\ell_{T_i})_{i\in\NN}$ converge vers $\ell$ dans ${\RR_+}^\C$, alors
il existe un arbre (réel) $T$, muni d'une action isométrique de $\pi_1S$
à petits stabilisateurs d'ar\^ete, tel que $\ell=\ell_T$.

\blemm \label{lem:disttransdiff}
Si $T$ est un arbre (réel) $T$, muni d'une action isométrique
de $\pi_1S$ à petits stabilisateurs d'ar\^ete, et si $\sigma$
appartient à $\ts$, alors $\ell_T\neq \ell_\sigma$.  
\elemm

\dem Il existe (voir par exemple \cite{CM,Pau3}) deux éléments $g$ et
$h$ de $\pi_1S$ qui sont hyperboliques dans $T$ et d'axes de
translation disjoints. Donc $g$ et $h$ sont des éléments non triviaux,
non paraboliques dans $\pi_1S$, et n'appartiennent pas à un même
groupe cyclique. De plus (voir par exemple \cite{CM,Pau3}), pour tous
$n,m$ dans $\ZZ-\{0\}$, la quantité
$\ell_T(g^nh^m)-\ell_T(g^n)-\ell_T(g^m)$ est égale à la distance entre
les axes de translation de $g$ et de $h$ dans $T$, donc est
strictement positive et indépendante de $n$ et $m$.

Notons $A_g$ et $A_h$ les axes de translations de $g$ et $h$
respectivement dans $(\wt S,\wt \sigma)$. S'ils sont disjoints dans
$\wt S$, alors ils n'ont pas de point commun à l'infini par discrétude
de $\rho_\sigma(\pi_1S)$, et d'après \cite[Theo.~7.38.3]{Bea}, pour
tous $n,m$ dans $\ZZ-\{0\}$ avec $|n|$ et $|m|$ assez grands, {\small
$$\cosh \frac{1}{2}\ell_\sigma(g^nh^m)=\left|\;\cosh d(A_g,A_h)\;
\sinh \frac{1}{2}\ell_\sigma(g^n)\; 
\sinh \frac{1}{2}\ell_\sigma(h^m)+
\epsilon\cosh \frac{1}{2}\ell_\sigma(g^n)\; 
\cosh \frac{1}{2}\ell_\sigma(h^m)\;\right|\;,$$} 
 où $\epsilon$ vaut $+1$ ou $-1$ suivant que $g^n$ et $h^m$
translatent dans la m\^eme direction ou pas. Si $A_g$ et $A_h$ se
rencontrent dans $\wt S$, comme $g$ et $h$ ne sont pas dans un m\^eme
groupe cyclique, leurs axes de translation se coupent avec un angle
$\theta$ tel que $0<\theta<\pi$, et d'après \cite[Theo.~7.38.6]{Bea},
{\small
$$\cosh \frac{1}{2}\ell_\sigma(g^nh^m)=\epsilon\cos \theta\;
\sinh \frac{1}{2}\ell_\sigma(g^n)\; \sinh \frac{1}{2}\ell_\sigma(h^m)+
\cosh \frac{1}{2}\ell_\sigma(g^n)\; \cosh
\frac{1}{2}\ell_\sigma(h^m)\;,$$}
 où $\epsilon$ vaut $+1$ ou $-1$ suivant que $\theta$ mesure
l'angle entre les points attractifs de $g^n$ et de $h^m$ ou l'angle
entre le point attractif de $g^n$ et le point répulsif de $h^m$. Si
$\ell_T= \ell_\sigma$, ces deux formules contredisent le fait que
$$\ell_\sigma(g^nh^m)-\ell_\sigma(g^n)-\ell_\sigma(g^m)$$ 
soit indépendants de $n,m$ dans $\ZZ-\{0\}$.  
\cqfd

\medskip 
Puisque $\PP(\RR_+\mbox{}^{\C})$ est métrisable, ce lemme et
l'alinéa le précédant montrent que l'image de $\pi\circ\ell$ est
ouverte dans son adhérence, car ils impliquent que pour toute suite
convergente dans la frontière de l'image, sa limite n'est pas dans
l'image.

Le théorème \ref{theo:benouaisquoi} en découle. \cqfd

\medskip
\begin{center}
+ --------------------------- +
\end{center}

\medskip 
Montrons maintenant comment construire la compactification de
Thurston de l'espace de Teichmüller directement par la topologie de
Gromov équivariante.

Soit ${\rm Hyp}(S)$ l'ensemble des classes d'isom\'etrie
\'equivariante de couples $(X,\alpha)$, o\`u $X$ est une vari\'et\'e
riemannienne compl\`ete simplement connexe de dimension $2$ \`a
courbure sectionnelle strictement n\'egative constante, et $\alpha$
une action propre et libre de $\pi_1 S$ telle que tout élément
parabolique de $\pi_1 S$ agisse par une isométrie parabolique sur $X$.
Munissons ${\rm Hyp}(S)$ de la topologie de Gromov équivariante.

Soit ${\rm Arb}(S)$ l'ensemble des classes d'isom\'etrie
\'equivariante de couples $(X,\alpha)$, o\`u $X$ est un arbre (réel)
et $\alpha$ une action isométrique de $\pi_1 S$ sur $X$, à petits
stabilisateurs d'arête (donc en particulier sans point fixe global, et
telle que tout élément parabolique de $\pi_1 S$ agisse par une
isométrie elliptique sur $X$), et {\it minimale} (i.e.~sans sous-arbre
invariant propre non trivial). (Cet ensemble porte d'autres noms dans
la littérature, voir par exemple \cite{CM,Sha1,Sha2,Bes2,Mor,Pau4}).
Munissons ${\rm Arb}(S)$ de la topologie de Gromov équivariante.

Munissons l'ensemble somme disjointe ${\rm Hyp}(S)\sqcup{\rm Arb}(S)$
de la topologie de Gromov équivariante, qui induit sur chacun des
sous-ensembles ${\rm Hyp}(S)$ et ${\rm Arb}(S)$ leur topologie de
Gromov équivariante. L'action de ${\rm Mod}(S)$ sur cet ensemble,
définie par l'application $(f,(X,\alpha))\mapsto (X,\alpha\circ f_*)$,
est une action par homéomorphismes. Le groupe topologique $\RR_+^*$
agit continûment sur ${\rm Hyp}(S)\sqcup{\rm Arb}(S)$ par
l'application $(t,(X,\alpha))\mapsto(tX,\alpha)$, et les actions de
$\RR_+^*$ et de ${\rm Mod}(S)$ commutent. Notons $${\rm K}(S)=(\,{\rm
Hyp}(S)\sqcup{\rm Arb}(S)\,)\,/\;\RR_+^*$$ l'espace topologique
quotient, muni de l'action quotient de ${\rm Mod}(S)$. La projection
canonique $p:{\rm Hyp}(S)\sqcup{\rm Arb}(S)\ra {\rm K}(S)$ est
ouverte. Notons $\Theta:{\rm Hyp}(S)\sqcup{\rm
Arb}(S)\ra{\RR_+}^\C-\{0\}$ l'application définie par
$$(X,\alpha)\mapsto ([\ga]\mapsto \min\{d(x,\alpha(\ga,x))\;:x\in
X\})\;,$$ qui est équivariante pour les actions de ${\rm Mod}(S)$, et
$\overline\Theta:{\rm K}(S)\ra\PP({\RR_+}^\C)$ l'application obtenue
par passage au quotient, qui est aussi équivariante pour les actions
de ${\rm Mod}(S)$.

Notons $\iota:\ts\ra {\rm Hyp}(S)\sqcup{\rm Arb}(S)$ l'application
définie par $\sigma\mapsto (\HH^2_\RR,\rho_\sigma)$, et
$\overline\iota:\ts\ra {\rm K}(S)$ l'application quotient, qui est
équivariante pour les actions de ${\rm Mod}(S)$.

\btheo (F.~Paulin \cite{Pau1}) L'espace ${\rm K}(S)$ est métrisable
compact. L'application $\overline\iota:\ts\ra {\rm K}(S)$ est un
homéomorphisme sur son image, celle-ci étant ouverte et dense dans
${\rm K}(S)$. Donc $(\,\overline\iota,{\rm K}(S)\,)$ est une
compactification de $\ts$ naturelle pour le groupe modulaire ${\rm
  Mod}(S)$.

L'application $\overline\Theta:{\rm K}(S)\ra\PP({\RR_+}^\C)$ est un
homéomorphisme sur son image, telle que
$\overline\Theta\circ\overline\iota=\pi\circ\ell$, donc la
compactification naturelle ci-dessus est isomorphe à la
compactification de Thurston.  
\etheo

En utilisant les travaux d'A.~Parreau \cite{Par}, un résultat
semblable pour des espaces symétriques de type non compact est
possible (sauf peut-\^etre la densité de l'image).

\medskip
\dem D'après un théorème d'Urysohn (voir par exemple \cite{Dug} page
233), pour montrer qu'un espace topologique $Y$ est métrisable compact,
il suffit de montrer qu'il est séparé, à base dénombrable d'ouverts,
et {\it séquentiellement compact} (i.e.~que de toute suite
(d\'enombrable) dans $Y$, on peut extraire une sous-suite
convergente).

Munissons $\RR_+^*\times \ts$ de la topologie produit des topologies
usuelles, et de l'action par homéomorphismes évidente de $\RR_+^*$,
définie par $t\cdot(t',\sigma)=(tt',\sigma)$, ainsi que celle de ${\rm
  Mod}(S)$, définie par $f\cdot (t',\sigma)=(t',f\cdot\sigma)$. En
utilisant la proposition \ref{pro:coincidtopos}, il n'est pas
difficile de montrer que l'application de $\RR_+^*\times \ts$ dans
${\rm Hyp}(S)$, qui à $(\lambda,\sigma)$ associe $(\lambda\HH^2_\RR,
\rho_\sigma)$ est un homéomorphisme. Elle induit par passage au
quotient un homéomorphisme de $\ts$ sur $p({\rm Hyp}(S))$ équivariant
pour les actions de ${\rm Mod}(S)$. En particulier, $p({\rm Hyp}(S))$
est séparable.

Il est montré dans \cite{Pau4} que l'application de ${\rm Arb}(S)$
dans ${\RR_+}^\C-\{0\}$, qui à $T$ associe $\ell_T$, est un homéomorphisme
sur son image. Elle commute avec les actions de $\RR_+^*$ et de ${\rm
  Mod}(S)$. (En fait, dans cette référence \cite{Pau4}, on ne
demandait pas que les distances de translation des éléments
paraboliques de $\pi_1S$ soient nulles, et on avait $\pi_1S$ à la
place de $\C$, et on travaillait avec des groupes plus généraux que
$\pi_1S$, mais le résultat énoncé ci-dessus découle de ce cas
général.)  Donc l'application $T\mapsto\ell_T$ induit un
homéomorphisme sur son image de $p({\rm Arb}(S))$ dans
$\PP({\RR_+}^\C)$, et son image est métrisable.  Par une preuve de
même schéma que pour le théorème \ref{theo:dichotomiequoi}, il est
facile de montrer que $p({\rm Arb}(S))$ est séquentiellement compact :
l'analogue du lemme \ref{lem:bestvina} est obtenu par le fait qu'une
action à petits stabilisateurs d'ar\^ete de $\pi_1S$ sur un arbre
(réel) ne fixe pas de point à l'infini de l'arbre (voir par exemple
\cite{CM,Pau3}); et l'analogue de la première contradiction dans la
preuve du lemme \ref{lem:limptstabar} est obtenue par le fait que si
$\eta$ est assez petit, alors $[\alpha,\beta]$ et $[\alpha,\beta^2]$
fixent un segment non trivial dans un arbre approchant. Donc $p({\rm
  Arb}(S))$ est métrisable compact (et en particulier séparable).
(Voir \cite{Pau2} pour une autre preuve; on retrouvait ainsi dans
\cite{Pau2} la compacité de l'image de $p({\rm Arb}(S))$ dans
$\PP({\RR_+}^\C)$, ce qui est un résultat de \cite[Theo.~5.3]{CM}.)

Tout arbre (réel) muni d'une action isométrique, sans point fixe
global et minimale, d'un groupe est réunion de ses axes de
translation (voir par exemple \cite{CM}), donc est séparable si le
groupe est dénombrable. Dans tout ensemble $\E$ d'espaces métriques
séparables munis d'une action isométrique d'un groupe dénombrable
$\Ga$, tout élément $X$ possède un système fondamental dénombrable de
voisinages pour la topologie de Gromov équivariante, puisqu'il suffit
de prendre les $V_{\epsilon,K,P}(X)$ pour $\epsilon$ dans $\QQ$ et $K$
une partie (finie) d'une partie dénombrable dense fixée de $X$.  Comme
$p$ est ouverte, tout élément de l'espace ${\rm K}(S)$ admet donc un
système fondamental dénombrable de voisinages. L'espace ${\rm K}(S)$
est de plus séparable, car $p({\rm Hyp}(S))$ et $p({\rm Arb}(S))$ le
sont.  Donc il est à base dénombrable d'ouverts.

Par un argument similaire à celui de la réciproque dans la preuve de
la proposition \ref{pro:coincidtopos}, l'application $\Theta$ est
continue.  Donc $\overline\Theta$ est continue. Elle est injective,
car injective en restriction à $p({\rm Hyp}(S))$ et $p({\rm Arb}(S))$,
et par le lemme \ref{lem:disttransdiff}. Comme $\overline\Theta$ est
continue, injective à valeurs dans un espace séparé, l'espace ${\rm
  K}(S)$ est séparé. (Voir aussi \cite[Chap.~IV.2]{Pau1} pour une
preuve directe de la séparation de la topologie de Gromov équivariante
sur l'ensemble ${\rm Hyp}(S)\sqcup{\rm Arb}(S)$.)

Remarquons que si une suite d'actions isométriques de $\pi_1S$ sur des
espaces métriques converge, pour la topologie de Gromov équivariante,
vers un arbre (réel) muni d'une action sans point fixe global, elle
converge aussi pour la topologie de Gromov équivariante vers tout
sous-arbre invariant, et en particulier vers son unique sous-arbre
minimal (qui est la réunion de ses axes de translations, voir par
exemple \cite{CM,Pau3}).

L'espace ${\rm K}(S)$ est séquentiellement compact, car $p({\rm
Arb}(S))$ l'est, et par le théorème \ref{theo:dichotomiequoi} et la
remarque ci-dessus.

Par le théorème d'Urysohn suscité, ${\rm K}(S)$ est donc compact.
L'application $\overline\Theta$, qui est continue injective d'un
espace compact dans un espace séparé, est donc un homéomorphisme sur
son image. Comme $p({\rm Arb}(S))$ est compact dans ${\rm K}(S)$
séparé (ou parce qu'une limite d'une suite d'actions isométriques d'un
groupe fixé sur des arbres (réels) est encore une action isométrique
de ce groupe sur un arbre (réel)), le sous-espace $p({\rm Arb}(S))$
est fermé, et donc $p({\rm Hyp}(S))$ est ouvert. Donc $\overline\iota$
est un homéomorphisme sur son image. Par un théorème de Skora
\cite{Sko} (voir aussi \cite{Ota}, ce résultat n'était pas disponible
au moment de l'écriture de \cite{Pau1}), l'image de $\overline\iota$
est dense dans ${\rm K}(S)$. Par conséquent, le couple
$(\,\overline\iota, \,{\rm K}(S)\,)$ est bien une compactification de
$\ts$. La naturalité et la relation
$\overline\Theta\circ\overline\iota=\pi\circ\ell$ étant claires par
construction, le résultat en découle.  \cqfd

\bigskip Nous n'avons décrit dans ces notes que l'aspect topologique
de la compactification de Thurston de l'espace de Teichmüller.  Une
interprétation des points du bord comme des feuilletages
transversalement mesurés (ou des laminations géodésiques
transversalement mesurées), duale à l'interprétation par les arbres
réels, est cruciale pour de plus amples informations (voir par exemple
\cite{FLP,Bon3}).

\bigskip
{\small
\noindent 
\begin{tabular}{l}
Département de Mathématiques et Applications\\
UMR 8553 CNRS\\
École Normale Supérieure\\
45 rue d'Ulm\\
75230 PARIS Cedex 05, FRANCE.\\
{\it e-mail: Frederic.Paulin@ens.fr}
\end{tabular}
}
 
\end{document}